\documentclass[a4paper]{article}
\usepackage{amsmath,amsfonts,amssymb}
\usepackage{multicol,multirow,diagbox,bm}
\usepackage{graphicx,subfigure}
\usepackage{lineno,cases}
\allowdisplaybreaks 
\usepackage{pgf,tikz}
\usetikzlibrary{arrows,shapes,automata,positioning}
\usetikzlibrary{fit}
\usepackage{pst-node}
\usepackage{cite}
\usepackage{color,xcolor}
\renewcommand{\bar}{\overline}
\renewcommand{\top}{\intercal}

\newcommand{\R}{\mathbb{R}}

\usepackage{bm} 

\newtheorem{assj}{Assumption}
\newtheorem{lemj}{Lemma}
\newtheorem{thmj}{Theorem}

\newtheorem{defn}{Definition}
\newtheorem{remj}{Remark}

\newcommand{\pb}{\noindent\textbf{Proof.} }
\newcommand{\pe}{\hfill\rule{4pt}{8pt}}

\usepackage{cite}
\begin{document}
	\title{Nash Equilibrium Seeking Over Directed Graphs 
		\thanks{This work was partially supported by the National Natural Science Foundation of China under Grants 61973043, 62003239, and 61703368, Shanghai Sailing Program under Grant 20YF1453000, Shanghai Municipal Science and Technology Major Project No. 2021SHZDZX0100, and Shanghai Municipal Commission of Science and Technology Project No. 19511132101.}
	}
	\author{Yutao Tang, Peng Yi,  Yanqiong Zhang, and  Dawei Liu
		\thanks{Y. Tang is with the School of Artificial Intelligence, Beijing University of Posts and Telecommunications, Beijing 100876, China. P. Yi is with the Department of Control Science and Engineering, Tongji University, Shanghai, 200092, China and  Shanghai Institute of Intelligent Science and Technology, Tongji University, Shanghai, 200092, China. Y. Zhang is with the School of Automation, Hangzhou Dianzi University, Hangzhou, 310018, China. D. Liu is with the China Research and Development Academy of Machinery Equipment, Beijing, 100086, China. (E-mails: yttang@bupt.edu.cn, yipeng@tongji.edu.cn, yqzhang@hdu.edu.cn, bit\_wei123@sina.com)}
	}	
	\date{ }
	\maketitle

	{\noindent\bf Abstract}: In this paper, we aim to develop distributed continuous-time algorithms over directed graphs to seek the Nash equilibrium in a noncooperative game. Motivated by the recent consensus-based designs, we present a distributed algorithm with a proportional gain for weight-balanced directed graphs. By further embedding a distributed estimator of the left eigenvector associated with zero eigenvalue of the graph Laplacian,  we extend it to the case with arbitrary strongly connected directed graphs having possible unbalanced weights. In both cases, the Nash equilibrium is proven to be exactly reached with an exponential convergence rate. An example is given to illustrate the validity of the theoretical results.
	

\maketitle

\section{Introduction}

Nash equilibrium seeking in noncooperative games has attracted much attention due to its broad applications in multi-robot systems, smart grids, and sensor networks  \cite{fudenberg1991,basar2018handbook,maschler2020game}.  In such problems, each decision-maker/player has an individual payoff function depending upon all players' decisions and aims at reaching an equilibrium from which no player has incentive to deviate. Information that one player knows about others and  the information sharing structure among these players play a crucial role in resolving these problems.  In a classical full-information setting, each player has access information including its own objective function and the decisions taken by the other players in the game \cite{li1987distributed, basar1999dynamic, stankovic2011distributed}. As the decisions of all other agents can be not directly available due to the privacy concerns or communication cost, distributed designs only relying on each player's local information are of particular interest, and sustained efforts have been made to generalize the classical algorithms to this case via networked information sharing.

In multi-agent coordination literature, the information structure (or the information sharing topology) among agents is often described by graphs \cite{mesbahi2010graph}.  Following this terminology,  the Nash equilibrium seeking problem in the classical full-information setting involves a complete graph  where any two players can directly communicate with each other \cite{li1987distributed, basar1999dynamic,shamma2005dynamic,frihauf2011nash,scutari2014real}. A similar scenario is the case when this full-decision information is obtained via broadcasts from a global coordinator \cite{grammatico2017dynamic}. By contrast, distributed rules via local communication and computation do not require this impractical assumption on the information structure.  
	
To overcome the difficulty brought by the lack of full information, a typical approach is to leverage the consensus-based mechanism to share information via network diffusion \cite{olfati2007consensus, swenson2015empirical, lou2016nash,koshal2016distributed}. To be specific, each player maintains a local estimate vector of all players' decisions and updates this vector by an auxiliary consensus process with its neighbors. After that, the player can implement a best-response or gradient-play rule with the estimate of the joint decision. For example, the authors conducted an asynchronous gossip-based algorithm for finding a Nash equilibrium in \cite{salehisadaghiani2016distributed}. The two awake players will appoint their estimates as their average and then take a gradient step. Similar results have been delivered for general connected graphs by extending classical gradient-play dynamics \cite{gadjov2019passivity, ye2017distributed}.   Along this line, considerable progress has been made with different kinds of discrete-time or continuous-time Nash equilibrium seeking algorithm with or without coupled decision constraints even for nontrivial dynamic players \cite{liang2017distributed, zeng2019generalized, de2019distributed, yi2018distributed, romano2019dynamic, zhang2019distributed, deng2019distributed, tatarenko2020geometric}.  However, all these results except a few for special aggregative games heavily reply on the assumption that the underlying communication graph is undirected, which definitely narrows down the applications of these Nash equilibrium seeking algorithms.  
	
Based on the aforementioned observations, this paper is devoted to the solvability of the Nash equilibrium seeking problem for general noncooperative games over directed graphs. Moreover, we aim to obtain an exponential convergence rate. Note that the symmetry of information sharing structure plays a crucial role in both analysis and synthesis of existing Nash equilibrium seeking algorithms. However, the information structure will lose such symmetry over directed graphs, which certainly makes the considered problem more challenging. 
	
To solve this problem, we start from the recent work \cite{gadjov2019passivity}. In \cite{gadjov2019passivity}, the authors presented an augmented gradient-play dynamics and showed the dynamics converge to consensus on the Nash equilibrium exponentially fast under undirected and connected graphs. We will first develop a modified version of gradient-play algorithms for weight-balanced digraphs  by adding a proportional gain, and then extend it to the case with arbitrary strongly connected digraph by further embedding a distributed estimator of the left eigenvector associated with zero eigenvalue of the graph Laplacian. Under some similar assumptions on the cost functions as in \cite{gadjov2019passivity}, we show that the developed two algorithms can indeed recover the exponential convergence rate in both cases. Moreover, by adding such a free-chosen proportional gain parameter, we provide an alternative way to remove the extra graph coupling condition other than singular perturbation analysis as that in \cite{gadjov2019passivity}.  To the best knowledge of us, this is the first exponentially convergent continuous-time result to solve the Nash equilibrium seeking problem over general directed graphs. 

The remainder of this paper is organized as follows: Some preliminaries are presented in Sections \ref{sec:pre}. The problem formulation is given in Section \ref{sec:form}. Then, the main designs are detailed in Section \ref{sec:main}. Following that, an example is given to illustrate the effectiveness of our algorithms in Section \ref{sec:simu}. Finally, concluding remarks are given in Section \ref{sec:con}.

\section{Preliminaries}\label{sec:pre}

In this section, we present some preliminaries of convex analysis  \cite{ruszczynski2006nonlinear} and graph theory \cite{mesbahi2010graph} for the following analysis.

\subsection{Convex analysis}

Let $\R^n$ be the $n$-dimensional Euclidean space and $\R^{n\times m}$ be the set of all $n\times m$ matrices. ${\bf 1}_n$ (or ${\bf 0}_n$) represents an $n$-dimensional all-one (or all-zero) column vector and ${\bm 1}_{n\times m}$ (or ${\bm 0}_{n\times m}$) all-one (or all-zero) matrix. We may omit the subscript when it is self-evident. $\mbox{diag}(b_1,\,{\dots},\,b_n)$ represents an $n\times n$ diagonal matrix with diagonal elements $b_i$ with $i=1,\,{\dots},\,n$. $\mbox{col}(a_1,\,{\dots},\,a_n) = [a_1^\top,\,{\dots},\,a_n^\top]^\top$ for column vectors $a_i$ with $i=1,\,{\dots},\,n$.  For a vector $x$ and a matrix $A$, $\|x\|$ denotes the Euclidean norm and $\|A\|$ the spectral norm.

A function $f\colon \R^m \rightarrow \R $ is said to be convex if, for any $0\leq a \leq 1$ and $\zeta_1,\zeta_2 \in \R^m$, $f(a\zeta_1+(1-a)\zeta_2)\leq af(\zeta_1)+(1-a)f(\zeta_2)$.
It is said to be strictly convex if this inequality is strict whenever $\zeta_1 \neq \zeta_2$. A vector-valued function $\Phi \colon \R^m \rightarrow \R^m$ is said to be  $\omega$-strongly monotone, if for any $\zeta_1,\, \zeta_2 \in \R^m$,  $(\zeta_1-\zeta_2)^\top [\Phi(\zeta_1)-\Phi(\zeta_2)]\geq \omega \|\zeta_1-\zeta_2\|^2$.  Function $\Phi\colon \R^m \rightarrow \R^m$ is said to be  $\vartheta$-Lipschitz, if for any $ \zeta_1, \, \zeta_2 \in \R^m$, $\|\Phi(\zeta_1)-\Phi(\zeta_2)\|\leq \vartheta \|\zeta_1-\zeta_2\|$.
Apparently, the gradient of an $\omega$-strongly convex function is $\omega$-strongly monotone.

\subsection{Graph theory}

A weighted directed graph (digraph) is described by $\mathcal {G}=(\mathcal{N}, \mathcal {E}, \mathcal{A})$ with the node set $\mathcal{N}=\{1,\,{\dots},\,N\}$ and the edge set $\mathcal {E}$. $(i,\,j)\in \mathcal{E}$ denotes an edge from node $i$ to node $j$. The weighted adjacency matrix $\mathcal{A}=[a_{ij}]\in \mathbb{R}^{N\times N}$ is defined by $a_{ii}=0$ and $a_{ij}\geq 0$. Here $a_{ij}>0$ iff there is an edge $(j,\,i)$ in the digraph.  The neighbor set of node $i$ is defined as $\mathcal{N}_i=\{j\mid (j,\, i)\in \mathcal{E} \}$. A directed path is an alternating sequence $i_{1}e_{1}i_{2}e_{2}{\dots}e_{k-1}i_{k}$ of nodes $i_{l}$ and edges $e_{m}=(i_{m},i_{m+1}) \in\mathcal {E}$ for $l=1,2,{\dots},k$.   If there is a directed path between any two nodes, then the digraph is said to be strongly connected.   The in-degree and out-degree of node $i$ are defined by $d^{\mbox{in}}_i=\sum\nolimits_{j=1}^N a_{ij}$ and $d^{\mbox{out}}_i=\sum\nolimits_{j=1}^N a_{ji}$. A digraph is weight-balanced if $d^{\mbox{in}}_i=d^{\mbox{out}}_i$ holds for any $i=1,\,\dots,\,N$.  The Laplacian matrix of $\mathcal{G}$ is defined as $L\triangleq D^{\mbox{in}}-\mathcal{A}$ with $D^{\mbox{in}}=\mbox{diag}(d^{\mbox{in}}_1,\,\dots,\,d^{\mbox{in}}_N)$.  Note that $L{\bm 1}_N={\bm 0}_N$ for any digraph. When it is weight-balanced, we have ${\bm 1}_N^\top L={\bm 0}_N^\top$ and the matrix  $\mbox{Sym}(L)\triangleq \frac{L+L^\top}{2}$ is positive semidefinite. 

Consider a group of vectors $\{{\bf 1},\, {a_2},\, \dots,\,a_N\}$ with $a_i$ the $i$th standard basis vector of $\R^N$, i.e., all entries of $a_i$ are zero except the $i$-th, which is one. These vectors are verified to be linearly independent. We apply the Gram-Schmidt process to them and obtain a group of orthonormal vectors $\{\hat a_1,\,\dots,\,\hat a_N\}$. Let $M_1=\hat a_1\in \R^{N}$ and $M_2=[\hat a_2~\dots~\hat a_N]\in \R^{N\times (N-1)}$. It can be verified that  $M_1=\frac{1}{\sqrt{N}}{\bm 1}_N$,  $M_1^\top M_1=1$, $M_2^\top M_2=I_{N-1}$, $M_2^\top M_1={\bm 0}_{N-1}$,  and $M_1 M_1^\top + M_2 M_2^\top=I_{N}$. Then, for a  weight-balanced and strongly connected digraph, we can order the eigenvalues of $\mbox{Sym}(L)$ as $0=\lambda_1<\lambda_2\leq \dots\leq \lambda_N$ and further have $\lambda_2 I_{N-1}\leq M_2^\top \mbox{Sym}(L)M_2\leq \lambda_N I_{N}$.

\section{Problem formulation} \label{sec:form}

In this paper, we consider a multi-agent system consisting of $N$ agents labeled as $\mathcal{N}=\{1,\,\dots,\,N\}$. They play an $N$-player noncooperative game defined as follows: Agent $i$ is endowed with a continuously differentiable cost function ${J}_i(z_i,\,{\bm z}_{-i})$, where $z_i\in \R$ denotes the decision (or action) profile of agent $i$ and  ${\bm z}_{-i}\in \R^{N-1}$ denotes the decision profile of this multi-agent system except for agent $i$. In this game, each player seeks to minimize its own cost function $J_i$ by selecting a proper decision $z_i$. Here we adopt a unidimensional decision variable for the ease of presentation and multiple dimensional extensions can be made without any technical obstacles. 

The equilibrium point of this noncooperative game can be defined as in \cite{basar1999dynamic}.
\begin{defn}
	Consider the game $\mbox{G}=\{\mathcal{N},\,{J}_i,\,\R\}$.  A decision profile $z^*=\mbox{col}(z_1^*,\,\dots,\,z_N^*)$ is said to be a Nash equilibrium (NE) of the game $\mbox{G}$ if	$J_i(z_i^*,\,{\bm z}_{-i}^*)\leq J_i(z_i,\,{\bm z}_{-i}^*)$ for any $i\in \mathcal{N}$ and $z_i\in \R$.
\end{defn}

At a Nash equilibrium, no player can unilaterally decrease its cost by changing the decision on its own, and thus all agents tend to keep at this state. 
Denote $F(z)\triangleq \mbox{col}(\nabla_1 J_1(z_1,\,{\bm z}_{-1}),\,\dots,\,\nabla_N J_N(z_N,\,{\bm z}_{-N}))\in \R^N$ with $\nabla_i J_i(z_i,\,{\bm z}_{-i})\triangleq \frac{\partial }{\partial z_i}J_i(z_i,\,{\bm z}_{-i})\in \R$. Here $F$ is called the pseudogradient associated with $J_1,\,\dots,\,J_N$.  

To ensure the well-posedness of our problem, the following assumptions are made throughout the paper:
\begin{assj}\label{ass:convexity}
	For each $i\in \mathcal{N}$,  the function $J_i(z_i,\,{\bm z}_{-i})$ is twice continuously differentiable, strictly convex and radially unbounded in $z_i\in \R$ for any fixed ${\bm z}_{-i}\in \R^{N-1}$. 
\end{assj}
\begin{assj}\label{ass:monotone}
	The pseudogradient $F$ is $\underline l$-strongly monotone and $\bar {l}$-Lipschitz for two constants $\underline{l},\, \bar l>0$. 
\end{assj}

These assumptions have been used in \cite{gadjov2019passivity} and \cite{de2019distributed}. Under these assumptions, our game $\mbox{G}$ admits a unique Nash equilibrium $z^*$ which can be characterized by the equation $F(z^*)={\bm 0}$ according to Propositions 1.4.2 and 2.2.7 in \cite{facchinei2003finite}.  

In a full-information scenario when agents can have access to all the other agents' decisions, a typical gradient-play rule 
\begin{align*}
\dot{z}_i=-\nabla_i J_i(z_i,\, {\bm z}_{-i}),\quad i\in \mathcal{N}
\end{align*}
can be used to compute this Nash equilibrium $z^*$. In this paper, we are more interested in distributed designs and assume that each agent only knows the decisions of a subset of all agents during the phase of computation.  

For this purpose, a weighted digraph $\mathcal{G}=(\mathcal{N},\, \mathcal{E}, \,\mathcal{A})$ is used to describe the information sharing relationships among the agents with node set $\mathcal{N}$ and weight matrix $\mathcal{A}\in \R^{N\times N}$. If agent $i$ can get the information of agent $j$, then there is a directed edge from agent $j$ to agent $i$ in the graph with weight $a_{ij}>0$.  Note that agent $i$ may not have the full-information of ${\bm z}_{-i}$ except the case with a complete communication graph. Thus, we have a noncooperative game with incomplete partial information. This makes the classical gradient-play rule unimplementable.  

To tackle this issue, a consensus-based rule has been developed in \cite{gadjov2019passivity} and each agent is required to estimate all other agents' decisions and implement an augmented gradient-play dynamics:
\begin{align}\label{sys:generator-gad}
\dot{{\bf z}}^i&=-\sum\nolimits_{j=1}^N a_{ij}({{\bf z}}^i-{{\bf z}}^j)- R_i \nabla_i J_i({\bf z}^i)
\end{align}
where $R_i=\mbox{col}({\bm 0}_{i-1},\,1,\,{\bm 0}_{N-i})$ and ${{\bf z}}^i=\mbox{col}(z_1^i,\,\dots,\,z_N^{i})$. Here ${\bf z}^i\in \R^N$ represents agent $i$'s  estimate of all agents' decisions with $z_i^i=z_i$ and  ${\bf z}_{-i}^i=\mbox{col}(z_1^i,\,\dots,z^i_{i-1},\,z^i_{i+1},\,\dots,\, z^i_{N})$. Function $\nabla_i J_i({\bf z}^i)=\frac{\partial J_i}{\partial z_i^i}(z_i^i,\,{\bf z}_{-i}^i)$ is the partial gradient of agent $i$'s cost function evaluated at the local estimate ${\bf z}^i$. 

For convenience, we define an extended pseudogradient as ${\bm F}({\bf z})=\mbox{col}( \nabla_1 J_1({\bf z}^1),\,\dots,\, \nabla_N J_N({\bf z}^N))\in \R^N$ for this game $\mbox{G}$.  The following assumption on this extended pseudogradient $\bm F$ is made in \cite{gadjov2019passivity}:
\begin{assj}\label{ass:lip-extended}
	The extended pseudogradient $\bm F$ is $l_{F}$-Lipschitz with $l_{F}>0$.
\end{assj}

Let $l=\max\{\bar l,\, l_{F}\}$. According to Theorem 2 in \cite{gadjov2019passivity}, along the trajectory of system \eqref{sys:generator-gad}, ${\bm z}^i(t)$ will exponentially converge to $z^*$ as $t$ goes to $+\infty$ if  graph $\mathcal{G}$ is undirected and satisfies a strong coupling condition of the form: $\lambda_2> \frac{l^2}{\underline{l}}+l$. Note that the coupling condition might be violated in applications for a given game and undirected graph (since the scalars $\lambda_2$ and $\frac{l^2}{\underline{l}}+l$ are both fixed).  Although the authors in \cite{gadjov2019passivity} further relaxed this connectivity condition by some singular perturbation technique, the derived results are still limited to undirected graphs. 

In this paper, we assume that the information sharing graph is directed and satisfies the following condition: 
\begin{assj}\label{ass:graph}
	Digraph $\mathcal{G}$ is strongly connected. 
\end{assj}

The main goal of this paper is to exploit the basic idea of algorithm \eqref{sys:generator-gad} and develop effective distributed variants to solve this problem  for digraphs under Assumption \ref{ass:graph} including undirected connected  graphs as a special case. Since the information flow might be asymmetric in this case, the resultant equilibrium seeking problem is thus more challenging than the undirected case.

\section{Main result}\label{sec:main}
In this section, we first solve our Nash equilibrium seeking problem for the weight-balanced digraphs and then extend the derived results to general strongly connected ones with unbalanced weights.  

\subsection{Weight-balanced graph}

To begin with, we make the following extra assumption:
\begin{assj}\label{ass:balance}
	Digraph $\mathcal{G}$ is weight-balanced. 
\end{assj}

Motivated by algorithm \eqref{sys:generator-gad}, we propose a modified version of gradient-play rules for game $\mbox{G}$ as follows:
\begin{align}\label{sys:generator-balanced}
\dot{{\bf z}}^i&=-\alpha \sum\nolimits_{j=1}^N a_{ij}({{\bf z}}^i-{{\bf z}}^j)- R_i \nabla_i J_i({\bf z}^i)
\end{align}
where $R_i$, ${{\bf z}}^i$ are defined as above and $\alpha>0$ is a constant to be specified later. Putting it into a compact form, we have \begin{align}\label{sys:generator-balanced-compact}
\dot{{\bf z}}=-\alpha {\bf L}{\bf z}-R {\bm F}({\bf z})
\end{align}
where ${\bf z}=\mbox{col}({\bf z}^1,\,\dots,\,{\bf z}^N)$, $R=\mbox{diag}(R_1,\,\dots,\,R_N)$ and ${\bf L}=L\otimes I_N$ with the extended pseudogradient ${\bm F}({\bf z})$.

Different from algorithm \eqref{sys:generator-gad} and its singularly perturbed extension presented in \cite{gadjov2019passivity}, we add an extra parameter $\alpha$ to increase the gain of the proportional term ${\bf L} {\bf z}$. With this gain being large enough, the effectiveness of algorithm \eqref{sys:generator-balanced-compact} is shown as follows:

\begin{thmj}\label{thm:exp}
	Suppose Assumptions \ref{ass:convexity}--\ref{ass:balance} hold. Let $\alpha>\frac{1}{\lambda_2}(\frac{l^2}{\underline{l}}+l)$. Then, for any $i\in \mathcal{N}$, along the trajectory of system \eqref{sys:generator-balanced-compact}, ${\bm z}^i(t)$ exponentially converges to $z^*$ as $t$ goes to $+\infty$.
\end{thmj}
\pb We first show that at the equilibrium of system \eqref{sys:generator-balanced-compact}, $z_i$ indeed reaches the Nash equilibrium of game $\mbox{G}$. In fact, letting the righthand side of \eqref{sys:generator-balanced} be zero, we have $\alpha {\bf L}{\bf z}^*+ R {\bm F}({\bf z}^*)={\bm 0}$. Premultiplying  both sides by ${\bm 1}^\top_N \otimes I_N$ gives
\begin{align*}
{\bm 0}=\alpha ({\bm 1}^\top_N \otimes I_N)(L\otimes I_N){\bf z}^*+ ({\bm 1}^\top_N \otimes I_N) R {\bm F}({\bf z}^*)
\end{align*}
Using ${\bm 1}^\top_NL=0$ gives 
${\bm 0}=({\bm 1}^\top_N \otimes I_N) R {\bm F}({\bf z}^*)$. 
By the notation of $R$ and $\bm F$, we have ${\bm F}({\bf z}^*)=\bm 0$. This further implies that ${\bf L}{\bf z}^*=\bm 0$. Recalling the property of $L$ under Assumption \ref{ass:graph}, one can determine some $\theta \in \R^N$ such that ${\bf z}^*={\bm 1}\otimes \theta $. This means ${\bm F}({\bm 1}\otimes \theta )=\bm 0$ and thus $\nabla_i J_i(\theta_i,\,\theta_{-i})=0$, or equivalently, $F(\theta)=\bm 0$. That is, $\theta$ is the unique Nash equilibrium $z^*$ of $\mbox{G}$ and ${\bf z}^*={\bm 1}\otimes z^* $.

Next, we show the exponential stability of system \eqref{sys:generator-balanced-compact} at its equilibrium ${\bf z}^*={\bm 1}\otimes z^*$. For this purpose, we denote $\tilde  {\bf z}={\bf z}-{\bf z}^*$ and perform the coordinate transformation $\bar {\bf z}_1=(M_1^\top\otimes I_N)\tilde  {\bf z}$ and $\bar {\bf z}_2=(M_2^\top\otimes I_N)\tilde  {\bf z}$. It follows that
\begin{align*}
\dot{\bar{\bf z}}_1&=-(M_1^\top\otimes I_N)R \Delta\\
\dot{\bar {\bf z}}_2&=-\alpha [(M_2^\top L M_2)\otimes I_{N}]\bar {\bf z}_2-(M_2^\top\otimes I_N) R\Delta
\end{align*} 
where  $\Delta\triangleq {\bm F}({\bf z})- {\bm F}({\bf z}^*)$. 

Let $V(\bar {\bf z}_1,\,\bar {\bf z}_2)=\frac{1}{2}(\|\bar {\bf z}_1\|^2+\|\bar {\bf z}_2\|^2)$. Then, its time derivative along the trajectory of system \eqref{sys:generator-balanced-compact} satisfies that
\begin{align}\label{eq1:lem-generator}
\dot{V}&=-\bar {\bf z}_1^\top (M_1^\top\otimes I_N)R \Delta-\bar {\bf z}_2^\top (M_2^\top\otimes I_N) R\Delta\nonumber\\
&- \alpha \bar {\bf z}_2^\top  \{[M_2^\top  L  M_2]\otimes I_{N}\}\bar {\bf z}_2\nonumber\\
&=-\tilde {\bf z}^\top  R \Delta - \alpha \bar {\bf z}_2^\top  [(M_2^\top \mbox{Sym}(L) M_2)\otimes I_{N}]\bar {\bf z}_2\nonumber\\
&\leq -\alpha \lambda_2 \|\bar {\bf z}_2\|^2-\tilde {\bf z}^\top R \Delta
\end{align}
Since $\tilde {\bf z}=(M_1 \otimes I_N) \bar {\bf z}_1+(M_2 \otimes I_N) \bar {\bf z}_2\triangleq \tilde {\bf z}_1+\tilde {\bf z}_2$, we split $\tilde {\bf z}$ into two parts to estimate the above cross term and obtain that
\begin{align*}
-\tilde {\bf z}^\top R \Delta&=(\tilde {\bf z}_1+\tilde {\bf z}_2 )^\top R [{\bm F}(\tilde {\bf z}_1+\tilde {\bf z}_2+{\bf z}^*)- {\bm F}({\bf z}^*)]\\
&=- \tilde {\bf z}_1^\top R [{\bm F}(\tilde {\bf z}_1+\tilde {\bf z}_2+{\bf z}^*)-{\bm F}(\tilde {\bf z}_1+{\bf z}^*)]\\
&-\tilde {\bf z}_2^\top R [{\bm F}(\tilde {\bf z}_1+\tilde {\bf z}_2+{\bf z}^*)-{\bm F}(\tilde {\bf z}_1+{\bf z}^*)]\\
&-\tilde {\bf z}_1^\top R [{\bm F}(\tilde {\bf z}_1+{\bf z}^*)-{\bm F}({\bf z}^*)]\\  
&-\tilde {\bf z}_2^\top R [{\bm F}(\tilde {\bf z}_1+{\bf z}^*)-{\bm F}({\bf z}^*)]
\end{align*}
As we have ${\bm F}({\bm 1}_N\otimes y)=F(y)$ for any $y\in \R^N$, it follows by the strong monotonicity of $F$  that
\begin{align*}
&\tilde {\bf z}_1^\top R [{\bm F}(\tilde {\bf z}_1+{\bf z}^*)-{\bm F}({\bf z}^*)]\\
&\quad =\frac{\bar {\bf z}_1^\top}{\sqrt{N}} [{\bf F}({\bf 1}\otimes (\frac{\bar {\bf z}_1}{\sqrt{N}}+y^*))-{\bm F}({\bf 1}\otimes y^*)]\\
&\quad =\frac{\bar {\bf z}_1^\top}{\sqrt{N}} [F(y^*+\frac{{\bar {\bf z}}_1}{\sqrt{N}})-{F}(y^*)]\\
&\quad \geq \frac{\underline{l}}{N}\|{\bar {\bf z}}_1\|^2
\end{align*} 
where we use the identity $({\bf 1}^\top \otimes I_N) R=I_N$ and $\tilde {\bf z}_1^\top  R =\frac{\bar {\bf z}_1^\top }{\sqrt{N}}$. Note that $\|R\|=\|M_2\|=1$ by definition. This implies that $\|R^\top \tilde {\bf z}_2 \|\leq  \|\tilde {\bf z}_2\| =\|\bar {\bf z}_2\|$. Then, under Assumptions \ref{ass:monotone} and \ref{ass:lip-extended}, we have that
\begin{align}\label{eq2:lem-generator}
-\tilde {\bf z}^\top R \Delta&\leq \frac{2 l}{\sqrt{N}}\|\bar {\bf z}_1\|\|\bar {\bf z}_2\|+l\|\bar {\bf z}_2\|^2-\frac{\underline{l}}{N}\|\bar {\bf z}_1\|^2 \end{align}
Bringing inequalities \eqref{eq1:lem-generator} and \eqref{eq2:lem-generator} together gives
\begin{align}\label{eq3:lem-generator}
\dot{V}&\leq -\frac{\underline{l}}{N}\|\bar {\bf z}_1\|^2 -(\alpha \lambda_2- l)\|\bar {\bf z}_2\|^2+\frac{2 l}{\sqrt{N}}\|\bar {\bf z}_1\|\|\bar {\bf z}_2\| \nonumber \\ 
&=-\begin{bmatrix}
\|\bar {\bf z}_1\|& 
\|\bar {\bf z}_2\|
\end{bmatrix} 
A_\alpha \begin{bmatrix}
\|\bar {\bf z}_1\|\\
\|\bar {\bf z}_2\|
\end{bmatrix}
\end{align}
with $A_\alpha=\begin{bmatrix} \frac{\underline l}{N}&-\frac{l}{\sqrt{N}}\\ -\frac{l}{\sqrt{N}}&\alpha \lambda_2-l \end{bmatrix}$. 
When $\alpha>\frac{1}{\lambda_2}(\frac{l^2}{\underline{l}}+l)$, matrix $A_\alpha$ is positive definite. Thus, there exists a constant $\nu>0$ such that
\begin{align*}
\dot{V}&\leq -\nu V
\end{align*}
Recalling Theorem 4.10 in \cite{khalil2002nonlinear}, one can conclude the exponential convergence of ${\bf z}(t)$ to ${\bf z}^*$, which implies that ${\bf z}^i(t)$ converges to $z^*$ as $t$ goes to $+\infty$. The proof is thus complete.
\pe

\begin{remj}\label{rem:proportional-gain}
	Algorithm \eqref{sys:generator-balanced-compact} is a modified version of the gradient-play dynamics \eqref{sys:generator-gad} with an adjustable proportional control gain $\alpha$.  The criterion to choose $\alpha$ clearly presents a natural trade-off between the control efforts and graph algebraic connectivity. 
	By choosing a large enough $\alpha$, this theorem ensures the exponential convergence of all local estimates to the Nash equilibrium $z^*$ over weight-balanced digraphs and also provides an alternative way to remove the restrictive graph coupling condition presented in \cite{gadjov2019passivity}.  
\end{remj}

\subsection{Weight-unbalanced graph} 
In this subsection, we aim to extend the preceding design to  general strongly connected digraphs. 
In the following, we first modify \eqref{sys:generator-balanced-compact} to ensure its equilibrium as the Nash equilibrium of game $\mbox{G}$, and then implement it in a distributed manner by adding a graph imbalance compensator.

At first, we assume a left eigenvector of the Laplacian $L$ associated with the trivial eigenvalue is known and denoted by $\xi=\mbox{col}(\xi_1,\,\dots,\,\xi_N)$, i.e., $\xi^\top L={\bm 0}$. Without loss of generality, we assume $\xi^\top {\bf 1}=1$. Then, $\xi$ is componentwise positive by Theorem 4.16 in Chapter 6 of   \cite{berman1994nonnegative}. Here we use this vector $\xi$ to correct the graph imbalance in system \eqref{sys:generator-balanced} as follows:
\begin{align}\label{sys:generator-digraph-nominal}
\dot{{\bf z}}^i&=- \alpha \xi_i \sum\nolimits_{j=1}^N a_{ij}({{\bf z}}^i-{{\bf z}}^j)- R_i \nabla_i J_i({\bf z}^i)
\end{align}
Similar ideas can be found in \cite{lou2016nash} and \cite{hadjicostis2018distributed}. We put this system into a compact form
\begin{align}\label{sys:generator-digraph-nominal-vector-compact}
\dot{{\bf z}}=-\alpha {\bf L}_\Lambda{\bf z}-R {\bm F}({\bf z})
\end{align}
where $\Lambda=\mbox{diag}({\xi_1},\,\dots,\,{\xi_N})$ and ${\bf L}_\Lambda=\Lambda L \otimes I_N $. It can be easily verified  that $\Lambda L$ is the associated Laplacian of a new digraph $\mathcal{G}'$, which has the same connectivity topology as digraph $\mathcal{G}$ but with scaled weights, i.e.,  $a'_{ij}={\xi_i}{a_{ij}}$ for any $i,\,j\in \mathcal{N}$.  As this new digraph $\mathcal{G}'$ is naturally weight-balanced, we denote $\lambda'_2$ as the minimal positive eigenvalue of $\mbox{Sym}(\Lambda L)$. 

Here is an immediate consequence of Theorem \ref{thm:exp}.
\begin{lemj}\label{lem:digraph-nominal}
	Suppose Assumptions \ref{ass:convexity}--\ref{ass:graph} hold and let $\alpha>\frac{1}{\lambda_2'}(\frac{l^2}{\underline{l}}+l)$. Then, for any $ i\in \mathcal{N}$, along the trajectory of system \eqref{sys:generator-digraph-nominal-vector-compact}, ${\bf z}^i(t)$ exponentially converges to $z^*$ as $t$ goes to $+\infty$.
\end{lemj}
Note that the aforementioned vector $\xi$ is usually unknown to us for general digraphs. To implement our algorithm, we embed a distributed estimation rule of $\xi$ 
into system \eqref{sys:generator-digraph-nominal} as follows:
\begin{align}\label{sys:digraph-estimator}
\begin{split}
\dot{\bm \xi}^i&=-\sum\nolimits_{j=1}^N a_{ij}({{\bm  \xi}}^i-{{\bm \xi}}^j)
\end{split}
\end{align}
where ${\bm \xi}^i=\mbox{col}(\xi^i_1,\,\dots,\,\xi^i_N)$ with $\xi_i^i(0)=1$ and $\xi^{i}_j=0$ for any $j\neq i \in \mathcal{N}$.

Here the diffusion dynamics of ${\bm \xi}^i$ is proposed to estimate the eigenvector $\xi$ by $\mbox{col}(\xi_1^1,\,\dots,\,\xi_N^N)$.  The following lemma shows the effectiveness of \eqref{sys:digraph-estimator}.

\begin{lemj}\label{lem:digraph-estimation}
	Suppose Assumption \ref{ass:graph} holds. Then, along the trajectory of  system \eqref{sys:digraph-estimator}, $\xi^i_i(t)>0$ for any $t\geq 0$ and exponentially converges to $\xi_i$ as $t$ goes to $+\infty$. 
\end{lemj}
\pb Note that the matrix $-L$ is essentially nonnegative in the sense that $\kappa I-L$ is nonnegative for all sufficiently large constant $\kappa>0$. Under Assumption \ref{ass:graph}, matrix $-L$ is also irreducible. By Theorem 3.12 in Chapter 6 of \cite{berman1994nonnegative}, the matrix exponential $\mbox{exp}(-Lt)$ is componentwise positive for any $t\geq 0$.  As the evolution of ${\bm \xi}_i=\mbox{col}(\xi^1_i,\,\dots,\,\xi_i^N)$ is governed by $\dot{\bm \xi}_i=-L{\bm \xi}_i$ with initial condition ${\bm \xi}_i(0)=\mbox{col}({\bm 0},\,1,\,{\bm 0})$. Thus, ${\bm \xi}_i(t)=\mbox{exp}(-Lt){\bm \xi}_i(0)>0$ for any $t$. By further using Theorems 1 and 3 in \cite{olfati2007consensus}, we have that ${\xi}^i_i(t)$ exponentially converges to the value $\xi_i^*=\frac{\xi_i}{\sum_{j=1}^N \xi_j}$ for any $i\in \mathcal{N}$ as $t$ goes to $+\infty$. Since $\xi=\mbox{col}(\xi_1,\,\dots,\,\xi_N)$ is a left eigenvector of $L$ associated with eigenvalue $0$, one can easily verify that ${\xi^*}^\top L=0$. Under Assumption \ref{ass:graph}, $0$ is a simple eigenvalue of $L$. Then, there must be a constant $c\neq 0$ such that $\xi=c\xi^*$. Note that $\xi^\top {\bf 1}={\xi^*}^\top {\bf 1}=1$. One can conclude that $c=1$ and thus complete the proof.
\pe

The whole algorithm to seek the Nash equilibrium is presented as follows:
\begin{align}\label{sys:generator-digraph-final}
\begin{split}
\dot{{\bf z}}^i&=-\alpha \xi^i_i  \sum\nolimits_{j=1}^N a_{ij}({{\bf z}}^i-{{\bf z}}^j)- R_i \nabla_i J_i({\bf z}^i)\\
\dot{\bm \xi}^i&=-\sum\nolimits_{j=1}^N a_{ij}({{\bm  \xi}}^i-{{\bm \xi}}^j) 
\end{split} 
\end{align}
with $\xi_i^i(0)=1$ and $\xi^{i}_j=0$ for any $j\neq i \in \mathcal{N}$.

Bringing Lemmas \ref{lem:digraph-nominal} and \ref{lem:digraph-estimation} together, we provide the second main theorem of this paper.

\begin{thmj}\label{thm:digraph}
	Suppose Assumptions \ref{ass:convexity}--\ref{ass:graph} hold and let $\alpha>\frac{1}{\lambda_2'}(\frac{l^2}{\underline{l}}+l)$.  Then, for any $i\in \mathcal{N}$, along the trajectory of system \eqref{sys:generator-digraph-final}, ${\bf z}^i(t)$ exponentially converges to $z^*$ as $t$ goes to $+\infty$.
\end{thmj}
\pb First, we put the algorithm into a compact form:
\begin{align}\label{sys:generator-digraph-final-compact}
\begin{split}
\dot{{\bf z}}&=-\alpha {\bf L}_{\Lambda'}{\bf z}-R {\bm F}({\bf z})\\
\dot{\bm \xi}&=-\bf L {{\bm \xi}}
\end{split}
\end{align}
where ${\bf L}_{\Lambda'}=\Lambda' L \otimes I_N$ and $\Lambda'=\mbox{diag}({\xi_1^1},\,\dots,\, {\xi_N^N})$. From this, one can further find that the composite system consists of two subsystems in a cascaded form as follows:
\begin{align*}
\begin{split}
\dot{{\bf z}}&=-\alpha {\bf L}_{\Lambda}{\bf z}-R {\bm F}({\bf z})-\alpha [(\Delta_\Lambda L)\otimes I_N]{\bf z}\\
\dot{\bm \xi}&=-\bf L {{\bm \xi}}
\end{split}
\end{align*}
where ${\bf L}_{\Lambda}$ is defined as in \eqref{sys:generator-digraph-nominal-vector-compact} and $\Delta_\Lambda= \Lambda'-\Lambda$. Note that the term $\alpha [(\Delta_\Lambda L)\otimes I_N]{\bf z}$ can be upper bounded by $\gamma_p\mbox{exp}(-\beta_p t)\|\bf z\|$ for some positive constants $\gamma_p$ and $\beta_p$ according to Lemma \ref{lem:digraph-estimation}. By viewing $\alpha [(\Delta_\Lambda L)\otimes I_N]{\bf z}$ as a vanishing perturbation of the upper subsystem, the unperturbed ${\bf z}$-subsystem is globally exponentially stable at its equilibrium ${\bf z}^*={\bf 1}_N\otimes z^*$ by Lemma \ref{lem:digraph-nominal}. Recalling Corollary 9.1 in \cite{khalil2002nonlinear}, the whole algorithm \eqref{sys:generator-digraph-final-compact} is globally exponentially stable at its equilibrium.  This implies that along the trajectory of system \eqref{sys:generator-digraph-final-compact}, ${\bm z}^i(t)$ exponentially converges to $z^*$ as $t$ goes to $+\infty$.  The proof is thus complete.
\pe

\begin{remj}
	In contrast to the algorithm \eqref{sys:generator-balanced} with proportional gains in Theorem \ref{thm:exp}, this new rule \eqref{sys:generator-digraph-final} further includes a distributed left eigenvector estimator to compensate the imbalance of the graph Laplacian. Compared with those equilibrium seeking results in \cite{koshal2016distributed, ye2017distributed,gadjov2019passivity} for undirected graphs, the proportional control  and graph imbalance compensator together facilitate us to solve this problem for strongly connected digraphs including undirected graphs as its special case.
\end{remj}

\begin{figure}
	\centering
	\includegraphics[width=0.86\textwidth]{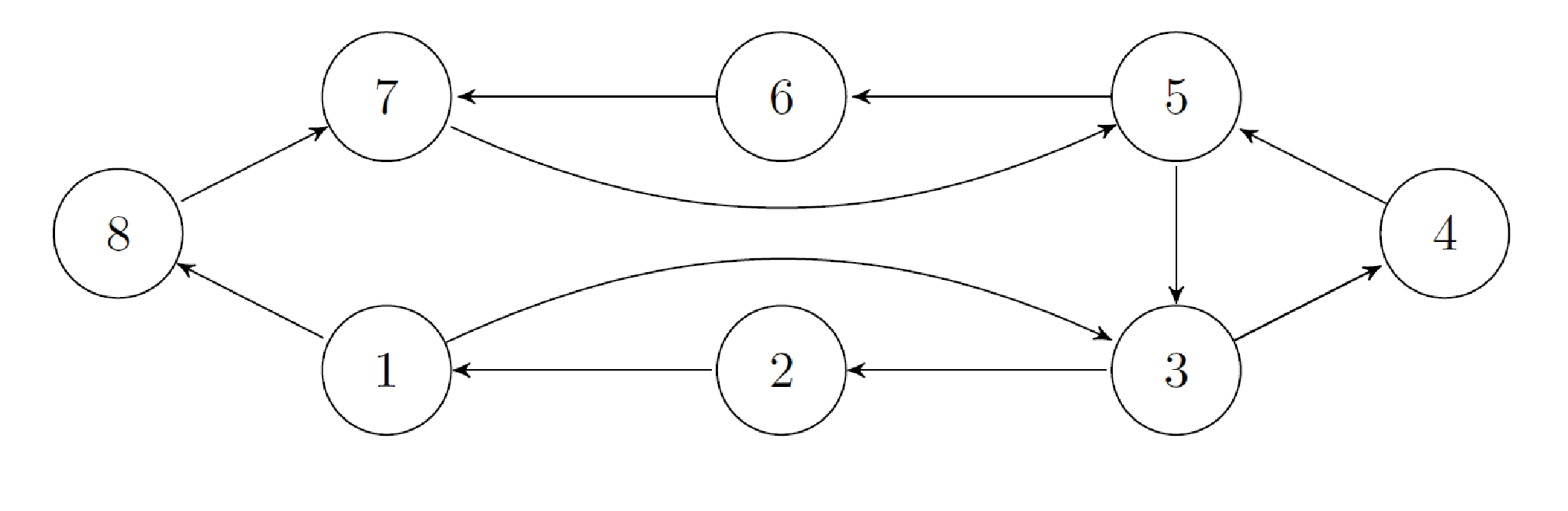}\\ 
	\caption{Digraph $\mathcal G$ in our example.}\label{fig:graph}
\end{figure}

\section{Simulation}\label{sec:simu}

In this section, we present an example to verify the effectiveness of our designs. 

Consider an eight-player noncooperative game. Each player has a pay-off function of the form $J_i(x_i,\,x_{-i})=c_i x_i -x_if(x)$ with $x=\mbox{col}(x_1,\,\dots,\,x_8)$ and $f(x)=D-\sum_{i=1}^8 x_i$ for a constant $D>0$. Suppose the communication topology among the agents is depicted by a digraph in Fig.~\ref{fig:graph} with all weights as one. The Nash equilibrium of this game can be analytically determined as $z^*=\mbox{col}(z_1,\,\dots,\,z_n)$ with $z^*_i=46-4*i$. 

Since the communication graph is directed and weight-unbalanced, the gradient-play algorithm developed in \cite{gadjov2019passivity} might fail to solve the problem. At the same time, Assumptions \ref{ass:convexity}--\ref{ass:graph} can be easily confirmed. Then, we can resort to Theorem \ref{thm:digraph} and use algorithm \eqref{sys:generator-digraph-final} to seek the Nash equilibrium in this eight-player noncooperative game.

\begin{figure}
	\centering
	\includegraphics[width=0.86\textwidth]{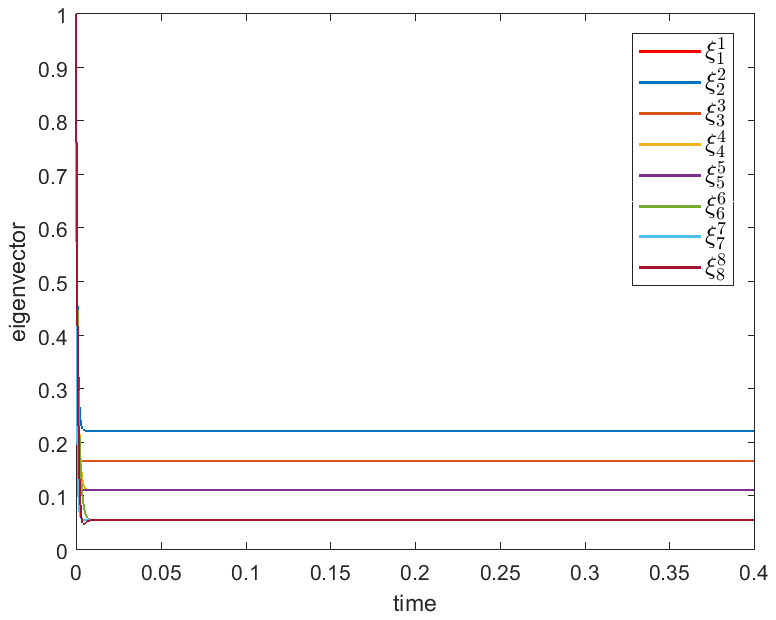}\\ 
	\caption{Profile of $\xi_i^i(t)$ in our example.}\label{fig:simu1}
\end{figure}

\begin{figure}
	\centering
	\includegraphics[width=0.86\textwidth]{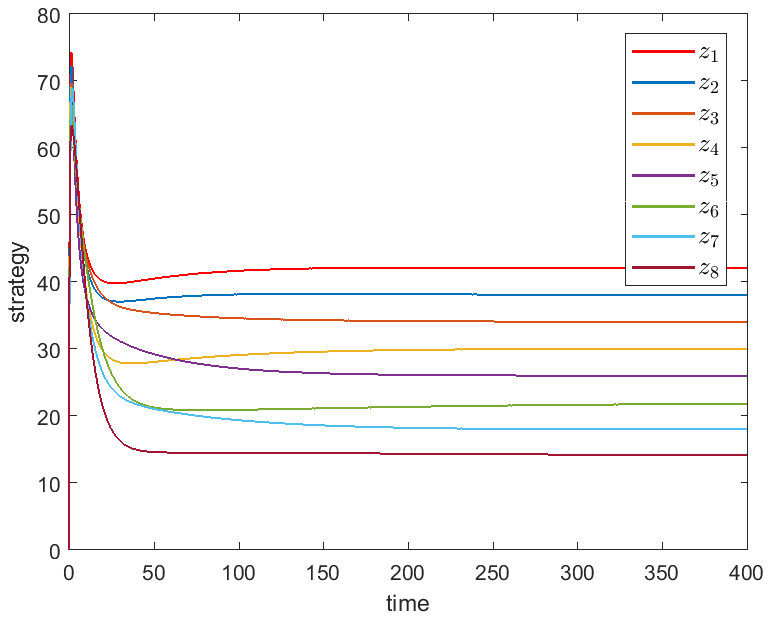}\\ 
	\caption{Profile of $z_i(t)$ in our example with $\alpha=2$.}\label{fig:simu2}
\end{figure}

\begin{figure}
	\centering
	\includegraphics[width=0.86\textwidth]{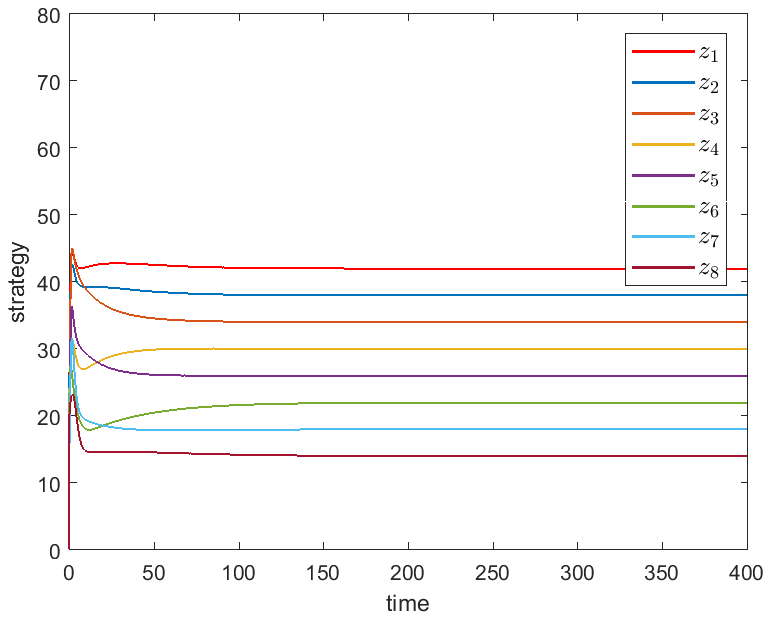}\\ 
	\caption{Profile of $z_i(t)$ in our example with $\alpha=10$.}\label{fig:simu3}
\end{figure}

\begin{figure}
	\centering
	\includegraphics[width=0.86\textwidth]{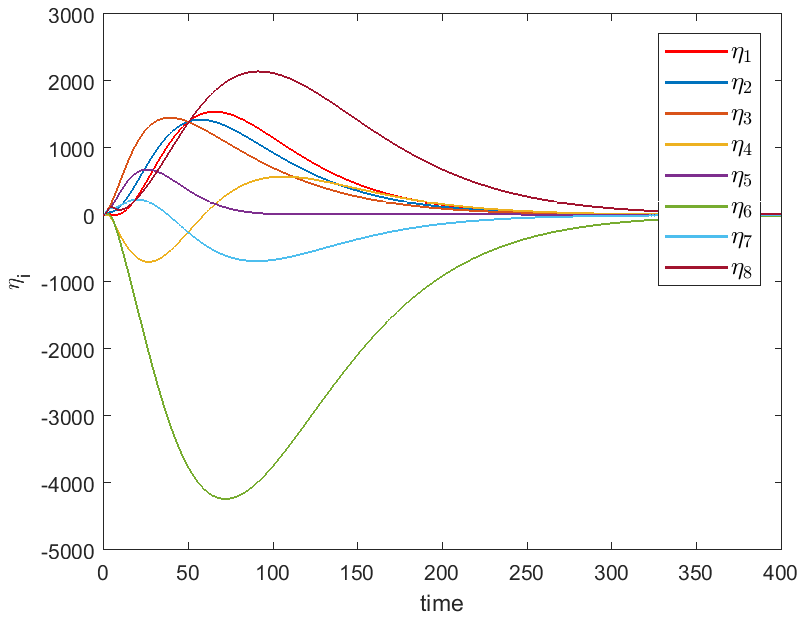}\\ 
	\caption{Profile of $\eta_i(t)$ in our example with $\alpha=10$.}\label{fig:simu4}
\end{figure}

For simulations, let $c_i=4i$ and $D=270$. We sequentially choose $\alpha=2$ and $\alpha=10$ for algorithm \eqref{sys:generator-digraph-final}. Since the righthand side of our algorithm is Lipschitz, we conduct the simulation via the forward Euler method with a small step size \cite{leveque2007finite}. The simulation results are shown in Figs.~\ref{fig:simu1}--\ref{fig:simu3}. From Fig.~\ref{fig:simu1}, one can find that the estimate $\xi(t)$ converges quickly to the left eigenvector of the graph Laplacian $\xi=\mbox{col}(4,\, 4,\, 3,\, 2,\, 2,\, 1,\, 1,\,1)/18$. At the same time, $\mbox{col}(z_1(t),\,\dots,\, z_8(t))$ approaches the Nash equilibrium $z^*$ of this game for different proportional parameters. Moreover, a larger proportional gain $\alpha$ is observed to imply a faster rate of convergence. We also show the profile of $\eta_i(t)\triangleq t^2(z_i(t)-z_i^*)$ in Fig.~\ref{fig:simu4} to confirm the exponential convergence rate when $\alpha=10$. These results verify the effectiveness of our designs in resolving the Nash equilibrium seeking problem over general strongly connected digraphs.

\section{Conclusion}\label{sec:con}
Nash equilibrium seeking problem over directed graphs has been discussed with consensus-based distributed rules. By selecting some proper proportional gains and embedding a distributed graph imbalance compensator,  the expected Nash equilibrium is shown to be reached exponentially fast over general strongly connected digraphs. In the future, we may use the adaptive high-gain techniques as in \cite{de2019distributed, tang2020optimal} to extend the results to fully distributed versions. Another interesting direction is to incorporate high-order agent dynamics and  nonsmooth cost functions.

\bibliographystyle{ieeetr}
\bibliography{game-exp}

\begin{thebibliography}{10}

\bibitem{fudenberg1991}
D.~Fudenberg and J.~Tirole, {\em Game Theory}.
\newblock Cambridge, USA: MIT Press, 1991.

\bibitem{basar2018handbook}
T.~Ba{\c{s}}ar and G.~Zaccour, {\em Handbook of Dynamic Game Theory}.
\newblock New York, USA: Springer, 2018.

\bibitem{maschler2020game}
M.~Maschler, S.~Zamir, and E.~Solan, {\em Game Theory}.
\newblock Cambridge, UK: Cambridge University Press, 2020.

\bibitem{li1987distributed}
S.~Li and T.~Ba{\c{s}}ar, ``Distributed algorithms for the computation of
  noncooperative equilibria,'' {\em Automatica}, vol.~23, no.~4, pp.~523--533,
  1987.

\bibitem{basar1999dynamic}
T.~Basar and G.~J. Olsder, {\em Dynamic Noncooperative Game Theory (2nd)}.
\newblock Philadelphia: SIAM, 1999.

\bibitem{stankovic2011distributed}
M.~S. Stankovic, K.~H. Johansson, and D.~M. Stipanovic, ``Distributed seeking
  of {Nash} equilibria with applications to mobile sensor networks,'' {\em IEEE
  Trans. Autom. Control.}, vol.~57, no.~4, pp.~904--919, 2011.

\bibitem{mesbahi2010graph}
M.~Mesbahi and M.~Egerstedt, {\em Graph Theoretic Methods in Multiagent
  Networks}.
\newblock Princeton, NJ, USA: Princeton University Press, 2010.

\bibitem{shamma2005dynamic}
J.~S. Shamma and G.~Arslan, ``Dynamic fictitious play, dynamic gradient play,
  and distributed convergence to {N}ash equilibria,'' {\em IEEE Trans. Autom.
  Control.}, vol.~50, no.~3, pp.~312--327, 2005.

\bibitem{frihauf2011nash}
P.~Frihauf, M.~Krstic, and T.~Basar, ``Nash equilibrium seeking in
  noncooperative games,'' {\em IEEE Trans. Autom. Control.}, vol.~57, no.~5,
  pp.~1192--1207, 2011.

\bibitem{scutari2014real}
G.~Scutari, F.~Facchinei, J.-S. Pang, and D.~P. Palomar, ``Real and complex
  monotone communication games,'' {\em IEEE Trans. Inf. Theory}, vol.~60,
  no.~7, pp.~4197--4231, 2014.

\bibitem{grammatico2017dynamic}
S.~Grammatico, ``Dynamic control of agents playing aggregative games with
  coupling constraints,'' {\em IEEE Trans. Autom. Control.}, vol.~62, no.~9,
  pp.~4537--4548, 2017.

\bibitem{olfati2007consensus}
R.~Olfati-Saber, J.~A. Fax, and R.~M. Murray, ``Consensus and cooperation in
  networked multi-agent systems,'' {\em Proc. IEEE}, vol.~95, no.~1,
  pp.~215--233, 2007.

\bibitem{swenson2015empirical}
B.~Swenson, S.~Kar, and J.~Xavier, ``Empirical centroid fictitious play: An
  approach for distributed learning in multi-agent games,'' {\em IEEE Trans.
  Signal Process.}, vol.~63, no.~15, pp.~3888--3901, 2015.

\bibitem{lou2016nash}
Y.~Lou, Y.~Hong, L.~Xie, G.~Shi, and K.~H. Johansson, ``Nash equilibrium
  computation in subnetwork zero-sum games with switching communications,''
  {\em IEEE Trans. Autom. Control.}, vol.~61, no.~10, pp.~2920--2935, 2016.

\bibitem{koshal2016distributed}
J.~Koshal, A.~Nedi{\'c}, and U.~V. Shanbhag, ``Distributed algorithms for
  aggregative games on graphs,'' {\em Oper. Res.}, vol.~64, no.~3,
  pp.~680--704, 2016.

\bibitem{salehisadaghiani2016distributed}
F.~Salehisadaghiani and L.~Pavel, ``Distributed {N}ash equilibrium seeking: A
  gossip-based algorithm,'' {\em Automatica}, vol.~72, pp.~209--216, 2016.

\bibitem{gadjov2019passivity}
D.~Gadjov and L.~Pavel, ``A passivity-based approach to {Nash} equilibrium
  seeking over networks,'' {\em IEEE Trans. Autom. Control.}, vol.~64, no.~3,
  pp.~1077--1092, 2019.

\bibitem{ye2017distributed}
M.~Ye and G.~Hu, ``Distributed {N}ash equilibrium seeking in multiagent games
  under switching communication topologies,'' {\em IEEE Trans. Cybern.},
  vol.~48, no.~11, pp.~3208--3217, 2017.

\bibitem{liang2017distributed}
S.~Liang, P.~Yi, and Y.~Hong, ``Distributed {N}ash equilibrium seeking for
  aggregative games with coupled constraints,'' {\em Automatica}, vol.~85,
  pp.~179--185, 2017.

\bibitem{zeng2019generalized}
X.~Zeng, J.~Chen, S.~Liang, and Y.~Hong, ``Generalized {Nash} equilibrium
  seeking strategy for distributed nonsmooth multi-cluster game,'' {\em
  Automatica}, vol.~103, pp.~20--26, 2019.

\bibitem{de2019distributed}
C.~De~Persis and S.~Grammatico, ``Distributed averaging integral {Nash}
  equilibrium seeking on networks,'' {\em Automatica}, vol.~110, p.~108548,
  2019.

\bibitem{yi2018distributed}
P.~Yi and L.~Pavel, ``Distributed generalized {N}ash equilibria computation of
  monotone games via double-layer preconditioned proximal-point algorithms,''
  {\em IEEE Trans. Control Netw. Syst.}, vol.~6, no.~1, pp.~299--311, 2018.

\bibitem{romano2019dynamic}
A.~Romano and L.~Pavel, ``Dynamic {NE} seeking for multi-integrator networked
  agents with disturbance rejection,'' {\em IEEE Trans. Control Netw. Syst.},
  vol.~7, no.~1, pp.~129--139, 2020.

\bibitem{zhang2019distributed}
Y.~Zhang, S.~Liang, X.~Wang, and H.~Ji, ``Distributed {Nash} equilibrium
  seeking for aggregative games with nonlinear dynamics under external
  disturbances,'' {\em IEEE Trans. Cybern.}, pp.~1--10, 2019.

\bibitem{deng2019distributed}
Z.~Deng and S.~Liang, ``Distributed algorithms for aggregative games of
  multiple heterogeneous {Euler}--{Lagrange} systems,'' {\em Automatica},
  vol.~99, pp.~246--252, 2019.

\bibitem{tatarenko2020geometric}
T.~Tatarenko, W.~Shi, and A.~Nedi{\'c}, ``Geometric convergence of gradient
  play algorithms for distributed {N}ash equilibrium seeking,'' {\em IEEE
  Trans. Autom. Control.}, vol.~66, no.~11, pp.~5342--5353, 2020.

\bibitem{ruszczynski2006nonlinear}
A.~Ruszczynski, {\em Nonlinear Optimization}.
\newblock Princeton: Princeton University Press, 2006.

\bibitem{facchinei2003finite}
F.~Facchinei and J.-S. Pang, {\em Finite-dimensional Variational Inequalities
  and Complementarity Problems}.
\newblock New York: Springer, 2003.

\bibitem{khalil2002nonlinear}
H.~K. Khalil, {\em Nonlinear Systems (3rd ed.)}.
\newblock New Jersey: Prentice Hall, 2002.

\bibitem{berman1994nonnegative}
A.~Berman and R.~J. Plemmons, {\em Nonnegative Matrices in the Mathematical
  Sciences}.
\newblock Philadelphia: SIAM, 1994.

\bibitem{hadjicostis2018distributed}
C.~N. Hadjicostis, A.~D. Dom{\'\i}nguez-Garc{\'\i}a, and T.~Charalambous,
  ``Distributed averaging and balancing in network systems: with applications
  to coordination and control.,'' {\em Found. Trends Syst. Control.}, vol.~5,
  no.~2-3, pp.~99--292, 2018.

\bibitem{leveque2007finite}
R.~J. LeVeque, {\em Finite Difference Methods for Ordinary and Partial
  Differential Equations}.
\newblock Philadelphia: SIAM, 2007.

\bibitem{tang2020optimal}
Y.~Tang and X.~Wang, ``Optimal output consensus for nonlinear multiagent
  systems with both static and dynamic uncertainties,'' {\em IEEE Trans. Autom.
  Control.}, vol.~66, no.~4, pp.~1733--1740, 2020.

\end{thebibliography}

\end{document}